\input amstex
\documentstyle{amsppt}
\magnification=1200

\catcode`\@=12


\pagewidth{4.8125 in}
\pageheight{7 in}
\hoffset=.2 in
\voffset=.3 in

\def\leaderfill{\leaders \hbox to 1em{\hss. \hss}\hfill}
\def\bs{\backslash}
\def\Lc{{\Cal L}}
\redefine\R{{{\bold R}}}
\redefine\C{{{\bold C}}}

\redefine\n{{{\bold n}}}

\redefine\L{\text{L}_\C^2}

\def\fg{{\frak g}}
\def\fz{{\frak z}}
\def\fv{{\frak v}}
\def\fw{{\frak w}}
\def\fa{{\frak a}}
\def\fg{{\frak g}}

\def\fh{{\frak h}}
\def\fr{{\frak r}}

\def\so{{\frak{so}}}
\def\ot{{\frac{1}{2}}}

\def\grad{{\mathop{grad}}}
\define\Lap{\Delta}

\def\dim{\mathop{dim}}
\def\ker{{\mathop{ker}}}
\def\tr{{\mathop{tr}}}

\def\Ric{{\mathop{Ric}}}
\def\Iso{{\mathop{Iso}}}
\def\bh{{\bar{H}}}
\def\fb{{\bar{f}}}
\def\pmaq{{M(j,c,Q)}}
\def\pmar{{M(j,c,r)}}
\def\pmapq{{M(j',c,Q)}}

\def\pn{{N}}
\def\al{{\alpha}}
\def\hc{{\Cal H}}
\def\bgj{{\bar{G}(j,c)}}
\def\bt{{\bar{\tau}}}
\def\bhj{{\bar{H}(j)}}
\def\bgjp{{\bar{G}(j',c)}}
\def\gj{{G(j,c)}}
\def\gjp{{G(j',c)}}
\def\ba{{\bar{A}}}
\def\th{{\tilde{H}}}
\def\bp{{\bar{\pi}}}

\newcount\theoremnumber
\def\grow{\advance \theoremnumber by 1}
\def\theorem#1\par{\grow
\noindent {\bf Theorem}\quad {\the\theoremnumber} .\qquad {\sl #1}\par}

\topmatter

\title Isospectral deformations of negatively curved Riemannian manifolds with boundary which are not locally
isometric \endtitle

\rightheadtext{Isospectral Riemannian metrics}

\author Carolyn S. Gordon and Zoltan I. Szabo
\endauthor

\affil Dartmouth College and Lehman College
\endaffil

\address 
\flushpar Carolyn S. Gordon:  Dartmouth College,  Hanover, New Hampshire, \ 03755;
\newline carolyn.s.gordon@dartmouth.edu
\flushpar Zoltan I. Szabo: Herbert H. Lehman College, Bronx, New York, \ 10468;
\newline zoltan@alpha.lehman.cuny.edu
\endaddress



\abstract To what extent does the eigenvalue spectrum of the Laplace-Beltrami operator on a compact Riemannian manifold determine the geometry of the
manifold?  We present a method for constructing isospectral manifolds with different local geometry, generalizing an earlier technique.  Examples include
continuous families of isospectral negatively curved manifolds with boundary as well as various pairs of isospectral manifolds.  The latter illustrate that
the spectrum does not determine whether a manifold with boundary has negative curvature, whether it has constant Ricci curvature, and whether it has parallel
curvature tensor, and the spectrum does not determine whether a closed manifold has constant scalar curvature.\endabstract

\thanks The first author is partially supported by a grant from the National Science Foundation.\endthanks 

\endtopmatter


\document

\heading Introduction \endheading

    A fundamental question in
spectral geometry is the extent to which the spectrum of the Laplacian on a Riemannian 
manifold determines the geometry of the manifold. The only way to identify specific
geometric invariants which are not spectrally determined is through explicit
constructions of isospectral manifolds, i.e., manifolds whose Laplacians, acting on
smooth functions, have the same eigenvalue spectrum.  In the case of manifolds with boundary, one may consider the spectrum of the Laplacian acting on
functions satisfying either Dirichlet or Neumann boundary conditions.  We will say that two manifolds with boundary are {\it isospectral} if they
are both Dirichlet and Neumann isospectral.

 All examples of isospectral manifolds constructed
prior to 1992 as well as many of the more recent examples are locally isometric; see, for example,
\cite{BGG}, \cite{Bu}, \cite{DG},
\cite{GWW}, \cite{GW2,3}, \cite{Gt1,2}, \cite{I}, \cite{M}, \cite{Su},
\cite{V} or the expository articles \cite{Be}, \cite{Br},  \cite{G3}, or \cite{GGt}. 
These examples reveal various global
invariants which are not spectrally determined such as the diameter and the fundamental group, but give no information concerning local
invariants such as curvature.  

In the past several years, examples of isospectral manifolds with different local geometry have
appeared.  The first such examples were pairs of manifolds with boundary constructed by the second author (preprint 1992).  These examples
together with examples constructed later of closed manifolds appeared in \cite{Sz}.  Among the latter examples are a pair of isospectral
closed manifolds, one of which is homogeneous and the other not.  The examples in \cite{Sz} were proven isospectral first by explicit computation
of the spectra and later by construction of an intertwining operator between the Laplacians.  The first author developed a technique for
constructing isospectral manifolds with different local geometry in \cite{G1,2}.  This technique was further developed in \cite{GW4}, resulting
in continuous families of isospectral manifolds with boundary having different Ricci curvature, and in \cite{GGSWW}, resulting in continuous
families of isospectral closed manifolds whose scalar curvature functions have different maxima.  D. Schueth
\cite{Sch} modified this technique to construct the first examples of isospectral simply-connected closed manifolds, in fact isospectral deformations of
simply-connected closed manifolds.  All these examples of isospectral manifolds with different local geometry are principal torus bundles with totally
geodesic fibers.

In this article, we give new examples of isospectral manifolds with different local geometry.  The manifolds are again principal torus bundles but the fibers are not totally geodesic. 
The new examples include:

\smallskip

(i)  Continuous isospectral deformations of negatively curved manifolds with boundary. The boundaries are also isospectral but in general non-isometric.  Examples
can be constructed in which the curvature is bounded above by any prescribed negative constant. The examples include both families of locally homogeneous manifolds and families of
locally inhomogeneous manifolds. These examples contrast with the result of Guillemin-Kazhdan
\cite{GuK} in dimension two, generalized by Croke and Sharafutdinov \cite{CS} to arbitrary dimensions, stating that negatively curved 
closed manifolds cannot be continuously isospectrally deformed.

(ii)  Pairs of isospectral manifolds with boundary one of which has negative sectional curvature and the other mixed curvature.

(iii)  A pair of isospectral manifolds with boundary one of which has constant Ricci curvature and the other variable Ricci curvature.

(iv)  Pairs of isospectral closed manifolds one of which has constant scalar curvature and the other variable scalar curvature.  In contrast, Patodi \cite{P}
showed that from the spectra of the Laplacian acting on functions, 1-forms, and 2-forms, one can tell whether the manifold has constant scalar curvature.  Our
examples show that the spectrum on functions alone does not determine this information.

(v)  Pairs of isospectral manifolds with boundary such that one manifold has parallel curvature tensor (it is a domain in a locally symmetric space) and the other
does not.  
 
\smallskip 

Examples (iii) and (iv) suggest the possibility that one may not be able to tell from the spectrum whether a closed manifold is Einstein, although the spectra
of the Laplacian on functions, 1-forms and 2-forms does determine whether the metric is Einstein \cite{P}.

It is very likely that the construction also allows a negatively curved manifold with boundary to be continuously isospectrally deformed to a manifold with curvature of both signs. 
The locally homogeneous isospectral manifolds $M_t$ with boundary constructed in (i) are domains with boundary in locally homogeneous manifolds whose universal
coverings are solvable Lie groups
$G_t$ with left-invariant metrics.  Each
$G_t$ is a semi-direct product $G_t=AH_t$ where $H_t$ is a nilpotent Lie group, the nilradical of the solvable group $G_t$, and $A\simeq \R$.  One may rescale the
metric on the factor $A$ arbitrarily (but consistently for the various $G_t$) without affecting the isospectrality of the domains $M_t$. For the particular class
of solvable Lie groups in our construction, Heintze
\cite{Hz} and, independently, Azencott and Wilson
\cite{AW2} proved that suitable scaling factors result in negative curvature.  In particular, for each of the
homogeneous solvmanifolds $G_t$, there exists a positive constant $\lambda_t$ such that rescaling the metric on $A$ by any constant less than $\lambda_t$ results
in negative curvature and any constant greater than $\lambda_t$ results in mixed curvature.  The isospectral manifolds with boundary cited in the second
example above were obtained by choosing a pair
$G_1$ and $G_2$ for which $\lambda_1$ and $\lambda_2$ are different.  (These pairs are not part of a continuous family.)  For these particular pairs, explicit
curvature computations are feasible and in fact were carried out in \cite{D}.  For the various continuous families $G_t$ which result in the first class of examples above,
explicit computations of the $\lambda_t$ seem intractible.  However, it is certainly likely that the $\lambda_t$ vary with $t$.  If this is the case, then one can
isospectrally deform a manifold of negative curvature to one of mixed curvature.
\medskip

\heading Section 1.  Technique for constructing isospectral manifolds with different local geometry.\endheading

\definition{1.1 Background and Notation}  Let $T$ be a torus, 
let $\pi:M\to N$ be a principal $T$-bundle and endow $M$ with a Riemannian metric so that the action of $T$ is by isometries.  
Give
$N$ the induced Riemannian metric so that
$\pi$ is a Riemannian submersion.  For
$a\in M$, we can decompose $T_a(M)$ into the vertical subspace (the tangent space to the fiber at $a$) and the horizontal subspace (the orthogonal complement to the vertical space).  For
$X\in T_a(M)$, write $X=X^v+X^h$ where $X^v$ is vertical and $X^h$ is horizontal.  Let $\nabla$ be the Levi-Civita connection on $M$.  The mean curvature of the fibers in $M$ is given
as follows:  Let $Z_1,\dots, Z_k$ be an orthonormal basis of left-invariant vector fields on $T$.  These define fundamental vector fields, which we denote by the same name, on $M$. 
These fundamental vector fields give a basis of the vertical space at each point.  For $a\in M$, the mean curvature at $a$  of the fiber through $a$ is given
by 
$$H_a=\Sigma_{i=1}^k (\nabla_{Z_i}Z_i)^h.$$
Since $T$ acts by isometries, we have $H_{z(a)}=z_*H_a$ for all $a\in M$ and $z\in T$.  Hence we may define a $\pi$-related vector field $\th$ on $N$.  We
will refer to $\th$ as the {\it mean curvature vector field of the submersion}.

Berard-Bergery and Bourguignon \cite{BB} gave a decomposition of the Laplacian $\Lap_M$ into vertical and horizontal components $\Lap_M =\Lap^v+\Lap^h$.  In the case of functions $f$ on
$M$ which are constant on the fibers of the submersion, so $f=\pi^*\fb$ for some function $\fb$ on $N$, then $\Lap^v(f)=0$ and  
$$\Lap_M(f)=\Lap^h(f)=\pi^*(\Lap_N(\fb)+\th(\fb)).\eqno{1.1}$$

If $N$ has non-trivial boundary, then $\partial M=\pi^{-1}(\partial N)$.  Since $\pi:M\to N$ is a Riemannian submersion, $\pi^*:C^\infty(N)\to C^\infty(M)$
maps functions on
$N$ satisfying Neumann boundary conditions to functions on $M$ satisfying Neumann boundary conditions.  Of course, it also maps functions satisfying Dirichlet conditions to functions
satisfying Dirichlet conditions.
\enddefinition

\remark{Remark}  In the situation of Notation 1.1, the space $\pi^*(C^\infty(N))$ of functions on $M$ which are constant on the fibers of the submersion is precisely the space of
$T$-invariant functions.  Since $T$ acts by isometries, it follows that $\pi^*(C^\infty(N))$ is invariant under $\Lap_M$.  The map $\pi^*:C^\infty(N)\to\pi^*(C^\infty(N))$ is a linear
isomorphism but it need not be unitary.  The operator $\Lap_N + \th$ on $N$ is not necessarily self-adjoint; however, it does have a discrete spectrum since it
is similar to the restriction of $\Lap_M$ to
$\pi^*(C^\infty(N))$.
\endremark

\proclaim{1.2. Theorem}  Let $T$ be a torus.  Suppose $M_1$ and $M_2$ are principal $T$-bundles endowed with Riemannian metrics so that the $T$-action is by
isometries.  For each subtorus
$K$ of
$T$ of codimension at most one, suppose that the quotient manifolds $K\bs M_1$ and
$K\bs M_2$, with the induced metrics, are isometric and that the isometry $\tau_K$ satisfies $\tau_{K*}(\th_K^{(1)})=\th_K^{(2)}$, where $\th_K^{(i)}$ is the mean curvature vector
field for the submersion $M_i\to K\bs M_i$.  (See Notation 1.1.)  Then
$M_1$ and
$M_2$ are isospectral.  If the manifolds have boundary, the conclusion is valid for both the Dirichlet and the Neumann spectrum.
\endproclaim

This theorem generalizes a method developed in \cite{G2} and \cite{GW4} in the special case that the fibers are totally geodesic.

\demo{Proof} Let $\Lap_i$ denote the Laplacian of $M_i$, and let $\L(M_i)$ denote the space of complex-valued square-integrable functions on $M_i$.  The torus
$T$ acts on
$\L(M_i)$,
$i=1,2$, and by a Fourier decomposition for this action, we have 
$$\L(M_i)=\Sigma _{\alpha\in \hat{T}}{\Cal H}_i^\alpha$$
where $\hat{T}$ consists of all characters on $T$, i.e., all homomorphisms from the group $T$ to the unit complex numbers, and 
$${\Cal H}_i^\alpha=\{f\in \L(M_i):z f=\alpha(z)f\text{ for all } z\in T\}.$$
Since the torus action on $M_i$ is by isometries, the Laplacian leaves each of the subspaces ${\Cal H}_i^\alpha$ invariant.  If $M_i$ has boundary, then we replace $\L(M_i)$ with the
subspace of functions satisfying either Neumann or Dirichlet boundary conditions (with the boundary conditions chosen consistently).  

Define an equivalence relation on $\hat{T}$ by $\alpha\equiv\beta$ if $\ker(\alpha)=\ker(\beta)$.  Let $[\alpha]$ denote the equivalence class of $\al$ and let $[\hat{T}]$ denote the
set of equivalence classes.  Setting 
$${\Cal H}_i^{[\alpha]}=\Sigma_{\beta\in [\al]}{\Cal H}_i^\beta ,$$
then 
$$\L(M_i)=\Sigma_{[\al]\in [\hat{T}]}{\Cal H}_i^{[\alpha]}.$$

For $\alpha =1$ the trivial character, we have $[1]=\{1\}$, and the space ${\Cal H}_i^1$ consists of those functions constant on the fibers of the submersion
$\pi_i: M_i\to T\bs M_i$.  By equation (1.1) and the remarks following 1.1, $\pi_i^*$ intertwines the restriction of $\Lap_i$ to $\hc_i^1$ with
$\bar{\Lap}_i+\th_T^{(i)}$ acting on
$\L(T\bs M_i)$, where $\bar{\Lap}_i$ denotes the Laplacian of $T\bs M_i$.  Moreover the isometry
$\tau_T$, whose existence is hypothesized in the theorem, gives a unitary isomorphism between
$\L(T\bs M_1)$ and
$\L(T\bs M_2)$ which intertwines
$\bar{\Lap}_1+\th_T^{(1)}$ with
$\bar{\Lap}_2+\th_T^{(2)}$.  Thus the  restrictions of the Laplacians $\Lap_i$ to the spaces ${\Cal H}_i^1$ are isospectral.

For non-trivial $\alpha\in\hat{T}$, the kernel  of $\alpha$ is a
subtorus $K$ of $T$ of codimension one.   The space
of all functions on $M_i$ constant on the fibers of the submersion $\pi_i:M_i\to K\bs M_i$ coincides with ${\Cal H}_i^{[\alpha]}\oplus {\Cal H}_i^1$.  We
can again use the hypothesis of the theorem to conclude that the restrictions of the Laplacians of $M_i$ to the subspaces ${\Cal
H}_i^{[\alpha]}\oplus {\Cal H}_i^1$, $i=1,2$ are isospectral.  Since we already know that the restrictions to the subspaces ${\Cal H}_i^1$ are isospectral, we
conclude that the restrictions to ${\Cal H}_i^{[\alpha]}$ are also isospectral, and the theorem follows.
\enddemo

We now introduce the specific classes of manifolds that will be our main objects of study.

\definition{1.3 Notation}

\noindent {\bf (i) Lie algebras $\fg(j)$, $\fh(j)$ and $\fr$.}  

 Starting with inner product spaces
$(\fv, \langle \ ,\ \rangle )$ and
$(\fz, \langle \ ,\ \rangle )$, a non-trivial linear map
$j:\fz\to\so(\fv)$, a one-dimensional vector space $\fa$ and a fixed choice of non-zero vector $A\in\fa$, we construct three Lie algebras $\fh(j)$, $\fg(j)$ and $\fr$ as
follows.  As vector spaces, we set 
$$\fh(j)=\fv+\fz$$
$$\fr=\fv +\fa$$
and
$$\fg(j)=\fv+\fz+\fa.$$  
Define a Lie bracket on $\fh(j)$ so that $\fz$ is central, $[\fv,\fv]\subset\fz$, and
$$\langle [X,Y],Z\rangle =\langle j(Z)X,Y\rangle $$
for $X,Y\in\fv$ and $Z\in\fz$.  This gives $\fh(j)$ the structure of a two-step nilpotent Lie algebra.  Next give $\fg(j)$ the unique bracket structure so that $\fh(j)$ is an ideal in
$\fg(j)$ and 
$$[A,X]=1/2X\,\,{\text{\rm and} }\,\,[A,Z]=Z$$
for $X\in\fv$ and $Z\in\fz$.  Define the Lie algebra structure on $\fr$ by declaring $\fv$ to be abelian and setting
$$[A,X]=1/2X$$
for $X\in\fv$.  With these structures, $\fg(j)$ is a solvable Lie algebra with nilradical $\fh(j)$, $\fz$ is an abelian ideal in $\fg(j)$, and $\fr$ is
isomorphic to the quotient
$\fz\bs\fg(j)$.  

\smallskip

\noindent{\bf (ii) Metric Lie algebras $\fg(j,c)$, $\fh(j)$ and $\fr(c)$.} 

 A  metric Lie algebra is a Lie algebra $\fg$ together with an inner product. If $G$ is any Lie group with Lie
algebra
$\fg$, then the inner product on $\fg$ defines a left-invariant Riemannian metric on $G$.  Two metric Lie algebras are said to be isomorphic if there exists a
Lie algebra isomorphism between them which is also an isometry with respect to their inner products.

In the notation of (i), the inner products on $\fv$ and $\fz$ define an inner product on $\fh(j)$ so that the
decomposition $\fh(j)=\fv+\fz$ is an orthogonal sum.  Given a real number $c>0$, define an inner product on $\fa$ by requiring that
$\|A\|=\frac{1}{c}$.  This inner product together with the inner products on $\fv$ and $\fz$ defines inner products on $\fg(j)$ and $\fr$ so
that the decompositions $\fg(j)=\fv+\fz+\fa$ and $\fr=\fa+\fv$ are orthogonal sums.  We will denote these metric Lie algebras by $\fg(j,c)$
and
$\fr(c)$.  As Lie algebras, $\fg(j,c)=\fg(j)$ and $\fr(c)=\fr$; only the metric depends on $c$.  We denote by
$G(j)$ and $R$ the simply-connected Lie groups with Lie algebras $\fg(j)$ and $\fr$, respectively, and by $G(j,c)$ and $R(c)$ the Lie groups $G(j)$ and $R$ endowed
with the left-invariant metrics defined by the inner product on the Lie algebra.  (We remark that the Riemannian manifold $R(1)$ is isometric to real hyperbolic space.) The
notation
$H(j)$ will be used both for the simply-connected Lie group with Lie algebra
$\fh(j)$ and for this Lie group endowed with the left-invariant metric corresponding to the inner product defined above on $\fh(j)$.

\smallskip

\noindent{\bf (iii)  Torus bundle $\pi_j:\bgj\to R(c).$}  

The Lie group exponential map
$\exp:\fg(j)\to G(j)$ is a diffeomorphism which restricts to a linear isomorphism from $\fz$ to an abelian normal subgroup $\exp(\fz)$ of $G(j)$, central in the nilradical
$H(j)$ but not central in $G(j)$.    
The quotient $\exp(\fz)\bs G(j)$ is isomorphic to the solvable Lie group $R$ and, with the metrics defined in $(ii)$, the projection $\pi_j:\gj\to R(c)$ is a
Riemannian submersion for each choice of $c$.  The metric induced on each fiber by the metric on $\gj$ is Euclidean.   

 Let $\Lc$ be a lattice of full rank in $\fz$ and let $T$ be the torus $\Lc\bs\fz$.  We identify
$\Lc$ with $\exp(\Lc)$, a discrete subgroup of $G(j)$ (central in $H(j)$ but not even normal in $G(j)$).  Let
$\bgj$ and
$\bhj$ be the quotients
$\Lc\bs\gj$ and
$\Lc\bs H(j)$, respectively, with the induced Riemannian metrics.  Since
$\Lc$ is central in
$H(j)$, the manifold $\bhj$ is a Lie group covered by $H(j)$ and its metric is homogeneous.  This is not the case for $\bgj$, although the metric is locally
homogeneous.  Moreover, the  torus $T:=\Lc\bs \exp(\fz)$ acts freely on
$\bgj$ by isometries and each orbit, with the induced metric, is a flat torus.  This action gives $\bgj$ the structure of a principal $T$ bundle with base
$R(c)$.  The projection $\pi_j:\bgj\to R(c)$ is a Riemannian submersion.  The restriction $\bp_j$ of $\pi_j$ to $\bh(j)$ is a Riemannian submersion
$\bp_j:\bh(j)\to\fv$ to the Euclidean space
$\fv$.

 Each left-invariant vector field
$X$ on
$G(j)$ descends to a vector field on
$\bgj$ which we will also denote by
$X$.  We will abuse language and refer to $X$ as an {\it invariant} vector field on $\bgj$ even though $\bgj$ is not homogeneous.
\smallskip

\noindent {\bf (iv) Riemannian submanifolds $Q(c)$ of $R(c)$, $M(j,c,Q)$ of $\bgj$, and $N_r(j)$ of $H(j)$.} 

Let $m=\dim(\fv)$. The  orthogonal group $O(m,\R)\simeq O(\fv)$ acts by
orthogonal automorphisms on the metric Lie algebra
$\fr(c)$ leaving
$\fa$ pointwise fixed.  The corresponding automorphisms of the Lie group $R$ are isometries with respect to all the metric structures $R(c)$.  Given any compact submanifold $Q$ of $R$
(with or without boundary) which is invariant under the action of
$O(m)$, let $Q(c)$ denote $Q$ endowed with the metric induced from that of $R(c)$.  Set 
$$M(j,c, Q)=\pi_j^{-1}(Q)$$ with the metric induced from $\bgj$.   
Then $M(j,c,Q)$ is a submanifold of $\bgj$ and $\pi_j:M(j,c,Q)\to Q(c)$ is a Riemannian submersion and a principal $T$ bundle.  The action of $T$ is by
isometries.

Let $S_r$, denote the geodesic sphere of radius $r$ in the Euclidean space $\fv$ and let
$N_r(j)=\bp_j^{-1}(S_r)$, where $\bp_j:\bh(j)\to\fv$ is the projection defined in (iii).  Thus we also have a Riemannian submersion $\bp_j:  N_r(j)\to
S_r$.  

\enddefinition

The manifolds $M(j,c,Q)$ will be  our main objects of study. The only cases of interest to us here will be when (a) $Q$ is a bounded domain in $R$ or (b) $Q$ is a submanifold of
codimension one in $R$ with or without boundary. The manifolds
$N_r(j)$ were studied in
\cite{GGSWW}.  We will use $N_r(j)$ 
 as an aid in understanding the manifolds $M(j,c, Q)$ in case (b).

In what follows, we will fix $\fv$, $\fz$, $\fa$, a choice of lattice $\Lc$ in $\fz$ and a choice of $Q$.  We will consider families of maps
$j_t:\fz\to\so(\fv)$ so that the manifolds in each family $\{M(j_t,c, Q)\}_t$ (with $c$ fixed) are isospectral.

\definition{1.4 Notation and Remarks }  We use the notation of 1.3, fixing a choice of $j$ and $\Lc$.  Subtori $K$ of the torus $T=\Lc\bs\fz$ correspond
to subspaces
$\fw$  of
$\fz$ spanned by lattice vectors in $\Lc$.  Given $K$ and thus $\fw$, let $\fz_K=\fz\ominus \fw$,  let $j_K=j_{|\fz_K}$, and let $\Lc_K$ be
the orthogonal projection of
$\Lc$ to
$\fz_K$.  Then $\Lc_K$ is a lattice of full rank in $\fz_K$.  Let $W$ be the connected normal
subgroup of $G(j,c)$ with Lie algebra $\fw$.  As in 1.3, use the data $(\fa, \fv, \fz_K, j_K, c, \Lc_K)$ to define a solvable Lie group $G(j_K,c)$, a quotient
manifold
$\bar{G}(j_K,c)=\Lc_K\bs G(j_K,c)$ and a submanifold $M(j_K,c,Q)$.  Then we have isometries:
$$G(j_K,c)\simeq W\bs G(j,c),$$ 
$$\bar{G}(j_K,c)\simeq K\bs\bgj$$ and
$$M(j_K,c,Q)\simeq K\bs M(j,c,Q)$$
where the manifolds on the right-hand-side of each equation have the quotient Riemannian metrics.  We will thus identify $\bar{G}(j_K,c)$ with $K\bs\bgj$ and
$M(j_K,c,Q)$ with $K\bs M(j,c,Q)$.

In 1.3, we fixed a choice of vector $A$ in $\fa$ and viewed $A$ both as a
left-invariant vector field on
$G(j,c)$ with norm $1/c$ and as an invariant vector field on
$\bgj$.  We continue to write $A$ for the analogous left-invariant vector field on $G(j_K,c)$
and the invariant vector field on $\bar{G}(j_K,c)$.  
\enddefinition

\proclaim{1.5. Lemma}  In the notation of 1.1, 1.3 and 1.4, let $K$ be a subtorus of $T$.  Then the
mean curvature vector field for the submersion $\bgj\to \bar{G}(j_K,c)$ is given by $c^2(\dim(K))A$, viewed as a vector field on
$\bar{G}(j_K,c)$.  The mean curvature vector field $\th_K$ for the submersion $\pmaq\to M(j_K,c,Q)$ is given by the orthogonal projection to $T(M(j_K,c,Q))$ of
$c^2(\dim(K))A_{|M(j_K,c,Q)}$.
\endproclaim

\demo{Proof}  Let $\nabla$ be the Levi-Civita connection on $\bgj$.  For invariant vector fields $U,V,W$ on $\bgj$ we have
$$2\langle \nabla_U(V),W\rangle =\langle [U,V],W\rangle +\langle [W,U],V\rangle +\langle [W,V],U\rangle .$$
Using this formula and the fact that $\|A\|\equiv 1/c$ with respect to the metric on $\bgj$, we see that for any unit vector $Z\in\fz$, we have 
$$\nabla_Z(Z)=c^2A.$$  Since $A$ is a horizontal vector field for the submersion $\bgj\to \bar{G}(j_K,c)$, the first statement follows from the formula for the
mean curvature given in 1.1.  The second statement is a consequence of the first.
\enddemo

\definition{1.6. Definition}  Let ${\frak v}$ and ${\frak z}$ be as above.  

(i) A pair $j,
j^{\prime}$ of linear maps from ${\frak z}$ to $\so({\frak v})$ will be called {\it
equivalent}, denoted
$j \simeq j^{\prime}$, if there exists orthogonal transformations $\alpha$ of $\fv$ and $\beta$ of $\fz$ such that
$$\al j(z)\al^{-1} = j^{\prime}(\beta(z))
$$
for all $z \in {\frak z}$.  

(ii) Let $\Lc$ be a lattice of full rank in $\fz$.  We will say that the pair $(j,\Lc)$ is equivalent to the pair $(j',\Lc ')$ if $j\simeq j'$ and if the map $\beta$
in definition (i) can be chosen so that $\beta(\Lc)=\Lc$.

(iii) The pair $j,
j^{\prime}$ will be called {\it isospectral}, denoted $j \sim
j^{\prime}$, if for each $z \in {\frak z}$, the eigenvalue spectra (with
multiplicities) of $j(z)$ and $j^{\prime}(z)$ coincide; i.e., for each $z\in\fz$, there exists an
orthogonal linear operator
$\alpha_z$ for which 
$$ {\al_z j(z) \al_z^{-1} = j^{\prime}(z)}.
$$
\enddefinition

\proclaim{1.7. Proposition}  We use the notation of 1.3, fixing $\fv$, $\fz$, $\fa$, $\Lc$ and $c\in\R$.  Let $j$ and $j^{\prime}$ be linear
injections from ${\frak z}$ to
$\so({\frak v})$.  In the notation of 1.6, the following are
equivalent:

\item{(a)} $j \simeq j'$;
\item{(b)}  $G(j,c)$ is isometric to $G(j',c)$;
\item{(c)} $\bgj$ is locally isometric to $\bgjp$

Moreover, if the pair $(j,\Lc)$ is equivalent to the pair $(j',\Lc ')$, then $\bgj$ is isometric to $\bgjp$.

\endproclaim

\demo{Proof} The local geometries of $G(j,c)$ and $\overline{G}(j,c)$ are identical. Thus $(c)$ is equivalent to saying that $G(j,c)$ is
locally isometric to
$G(j',c)$ which, by simple-connectivity, is equivalent to $(b)$.  

By Theorem 5.2 of \cite{GW1}, $\gj$ is isometric to $\gjp$ if and only if $\fg(j,c)$ and $\fg(j',c)$ are isomorphic as metric Lie algebras.  (See 1.3(ii).)
This condition is equivalent to the condition that $j \simeq j^{\prime}$.  The metric Lie algebra isomorphism $\tau:\fg(j,c)\to\fg(j',c)$ is given in this case
by 
$\tau(sA+X+Z)=sA+\alpha(X)+\beta(Z)$ for $s\in\R$,
$X\in\fv$ and $Z\in\fz$, where $\alpha$ and $\beta$ are given as in Definition 1.6(i).
The corresponding isomorphism
$\tilde{\tau}:G(j,c)\to G(j',c)$ is then an isometry. Thus (a) and (b) are equivalent.

Finally, if the pair $(j,\Lc)$ is equivalent to the pair $(j',\Lc ')$, then the isometry $\tilde{\tau}$ in (b) descends to an isometry $\tilde{\tau}:\bgj\to
\bgjp$.

\enddemo

\proclaim{1.8. Lemma}  We use the notation of 1.3(iv) and 1.6.  Suppose that the pair $(j,\Lc)$ is equivalent to the pair $(j',\Lc ')$.  If $Q$ is any
$O(m)$-invariant submanifold of
$R$, then $M(j,c,Q)$ is isometric to $M(j',c,Q)$.  
\endproclaim

\demo{Proof}  The isometry  $\tilde{\tau}:\bgj\to\bgjp$ defined in the proof of Proposition 1.7 is a bundle map. The induced isometry of $R(c)$ lies in $O(m)$ and is defined by the
orthogonal transformation $\alpha$ of $\fv$ in the notation of Definition 1.6.  In particular, this isometry carries $Q$ to $Q$, since $Q$ is $O(m)$-invariant.  Hence $\tilde{\tau}$
restricts to an isometry from $M(j,c,Q)$ to $M(j',c,Q)$.
\enddemo

We will prove a partial converse to Lemma 1.8 under a genericity condition in Proposition 1.10 below.  First, however, we address the question of isospectrality.

\proclaim{1.9. Theorem}  We use the notation of 1.3. fixing $\fv$, $\fz$, $\fa$, $\Lc$, $Q$ and $c$.  Let $j,j':\fz\to\so(\fz)$ be
isospectral linear maps as in Definition 1.6. 
Then $M(j,c, Q)$ is isospectral to $M(j',c, Q)$.  If these manifolds have boundary, then they are both Dirichlet and Neumann isospectral. 
Moreover, their boundaries, with the induced metrics, are also isospectral.
\endproclaim

\demo{Proof}  We apply Theorem 1.2.  Let $K$ be a subtorus of $T$ of co-dimension one and recall the Notation 1.4.  Since $\fz_K$ is one-dimensional, the
condition that $j$ be isospectral to
$j'$  implies that the pair $(j_K,\Lc_K)$ is equivalent to the pair $(j'_K,\Lc_K)$.  Consequently, by Lemma 1.8, there exists an isometry
$\sigma: M(j_K,c,Q)\to M(j'_K,c,Q)$. To see that
$\sigma_*$  carries the mean curvature vector field for the submersion
$M(j,c,Q)\to  M(j_K,c,Q)$ to that for the submersion $M(j',c,Q)\to  M(j'_K,c,Q)$, apply Lemma 1.5 together with the fact that $\sigma$ extends to an 
isomorphism 
$\sigma:\bar{G}(j_K,c,Q)\to \bar{G}(j'_K,c,Q)$ satisfying $\sigma_*(A)=A$ (as can be seen from the proofs of Proposition 1.7 and Lemma 1.8).  Thus the
manifolds $M(j,c, Q)$ and
$M(j',c, Q)$ satisfy the hypothesis of Theorem 1.2 for all co-dimension one subtori $K$ of
$T$.  A similar but easier argument shows that the hypothesis is also satisfied for $K=T$.  Theorem 1.2 therefore implies that
$M(j,c,Q)$ is isospectral to
$M(j',c,Q)$ (Dirichlet and Neumann isospectral if the boundaries are non-empty).   In the case where the manifolds have boundary, we have 
$\partial (M(j,c,Q))=M(j,c,\partial Q)$
and $\partial (M(j',c,Q))=M(j',c,\partial Q)$.  Since $\partial Q$ is an $O(m)$-invariant submanifold of $R$, 
the isospectrality of the boundaries follows from the main statement of the
theorem.
\enddemo

We now consider the converse of Lemma 1.8.

\proclaim{1.10 Proposition}  We use the notation of 1.3, fixing $\fv$, $\fz$, $\fa$, $\Lc$, $c$ and $Q$.  Assume that $Q$ is either 
a bounded domain in $R$ or a submanifold of
codimension one.  Suppose that
$j:\fz\to\so(\fv)$ satisfies the property that there are only finitely many orthogonal maps of $\fv$ which commute with all the transformations $j(Z)$, $Z\in\fz$.  Then if
$j':\fz\to\so(\fv)$ is any linear map for which $\pmaq$ is isometric to $\pmapq$, then $j \simeq j'$.
\endproclaim

The hypothesis on $j$ is generic.  When the hypothesis holds and when $Q$ is a domain with boundary in $R$, Proposition 1.10 may be applied
to $M(j,c,\partial Q)=\partial (M(j,c,Q))$ to conclude that $\partial (M(j,c,Q))$ is not isometric to $\partial (M(j',c,Q))$.  Thus Theorem 1.9 and
Proposition 1.10 together imply that $\pmaq$ and $\pmapq$ are isospectral manifolds with isospectral boundaries, but neither the manifolds
nor the boundaries are isometric.

  Before proving the proposition, we give a more explicit description of
$Q$.

\definition{1.11 Explicit construction of $R$ and $Q$} (a) The Lie group $R$ of Definition 1.3 is diffeomorphic to $\R^+\times\fv$.   The Lie group 
multiplication, viewed on
$\R^+\times \fv$, is given by 
$$(t,X)(t',X')=(tt',X+t^{1/2}X')$$
for $t,t'\in\R^+$ and $X,X'\in\fv$.  The left-invariant vector field $A$ is given by 
$$A(t,X)=\frac{d}{ds}_{|s=0}((t,X)(s,0))=t\frac{\partial}{\partial t}.$$  
Each $X\in\fv$ gives rise in a natural way to two different vector fields on $R$.  First, $X$ defines a left-invariant vector field, also denoted $X$. 
Secondly, ignoring the group structure and identifying the underlying manifold $R$ with the vector space $\R^+\times\fv$, then $X$
defines a directional derivative $D_X$ on $R$ given by
$D_X(g)(t,Y)=\frac{d}{ds}_{|s=0}g(t,Y+sX)$ for $g\in C^\infty(R)$.  The two vector fields are related by:
$$X(t,Y)=t^\ot D_X(t,Y).$$  

(b) To construct $Q$ explicitly, first suppose that $Q$ is a compact $O(m)$-invariant submanifold of $R$ of codimension one in $R$ with or without boundary. 
Then there exists an interval
$[t_1,t_2]$ in $\R^+$ and a  continuous function $f:[t_1,t_2]\to\R$ such that
$f$ is strictly positive and smooth on $(t_1,t_2)$ and 
$$Q=\{(a,X)\in R: a\in [t_1,t_2] \text{ and }\|X\|= f(a)\}.\eqno{(1.2)}$$
(Of course, when $Q$ is closed, $f$ must satisfy certain boundary conditions.  These conditions will never be used explicitly, however.)

Next suppose that $Q$ is an $O(m)$-invariant bounded domain in $R$.  Then there exists an interval $[t_1,t_2]$ in $\R^+$ and a function $f:[t_1,t_2]\to\R$ as above, satisfying the
additional condition that  
$f(t_1)=f(t_2)=0$, such that 
$$Q=\{(a,X): a\in [t_1,t_2] \text{ and }\|X\|\leq f(a)\}.$$
(Again, smoothness of $\partial Q$ imposes an additional boundary condition on $f$ which we will not need explicitly.)

\enddefinition

\demo{Proof of Proposition 1.10}  First suppose that $Q$ is a bounded domain in $R$ and thus $M(j,c,Q)$ and $M(j',c,Q)$ are bounded domains in $\bgj$ and
$\bgjp$, respectively.  Since the local geometry of $M(j,c,Q)$ is that of $\bgj$, the statement of the proposition is immediate from Proposition 1.7 in this
case.

We thus assume that $Q$ has codimension one in $R$. We first consider the mean curvature vector fields for the
submersions
$\pi_j:M(j,c,Q)\to Q$ and
$\pi_{j'}:M(j',c,Q)\to Q$. By Lemma 1.5, the two vector fields are identical and are given by $c^2k\ba$ where $k$ is the dimension of $\fz$ and
$\ba(m)$ is the orthogonal projection of $A(m)$ to
$T_m(Q)$ for
$m\in Q$.  
 
Define $f$ as in 1.11(b) and write $\phi(t,X)=\|X\|-f(t)$.  Using the local coordinate expressions for the left-invariant vector fields on $R$ given in
1.11(a), we obtain
$$\grad(\phi)(t,X)=-ct f'(t)(cA)+\frac{t^\ot X}{\|X\|}.$$
(We emphasize that on the left-hand-side of this equation, the element $X$ of $\fv$ is being used as one of the coordinates for a point in $R$, while on the
right-hand-side, it is being viewed as a left-invariant vector field.  Recall that both 
$cA$ and $\frac{X}{\|X\|}$ are unit vectors with respect to the Riemannian metric.  The gradient is computed with respect to the Riemannian metric.)
Thus a unit normal vector to $Q$ at $(t,X)$ is given by $-p(t)cA+q(t)\frac{X}{\|X\|}$ where 
$$p(t)=\frac{ct^\ot f'(t)}{\sqrt{1+c^2t(f'(t)^2}}\text{ and }
q(t)=\frac{1}{\sqrt{1+c^2t(f'(t)^2}}.$$  
The orthogonal projection $c\ba$ of $cA$ to $T_m(Q)$ at $m=(t,X)$ has length $q(t)$.  

Thus $\|\ba\|$ is constant on each ``slice'' $S(t)$ obtained by
holding $a=t$ with
$t$ in
$[t_1,t_2]$.  (See equation (1.2).)  Each such slice is a round sphere in
$\fv$ isometric to the sphere $S_{f(t)}$ defined in 1.3(iv).  The sphere may collapse to a single point when $t=t_1$ and/or $t=t_2$.  Letting
$\pi_j:M(j,c,Q)\to Q$ and
$\pi_{j'}:M(j',c,Q)\to Q$ be the bundle projections, then
$\pi_j^{-1}(S(t))$ and
$\pi_{j'}^{-1}(S(t))$ are submanifolds of $M(j,c,Q)$ and
$M(j',c,Q)$ isometric to $\pn_{f(t)}(j)$ and $\pn_{f(t)}(j')$, respectively. (See 1.3(iv) for the definition of $\pn_r(j)$.)

We first prove the proposition in the special case that $\tau:M(j,c,Q)\to  M(j',c,Q)$ is both an isometry and a bundle map with respect to the bundle
structures
$\pi_j:M(j,c,Q)\to Q$ and
$\pi_{j'}:M(j',c,Q)\to Q$.  The idea of the proof in this case will be to show that any  isometry between $M(j,c,Q)$ and $M(j',c,Q)$ induces
an  isometry between
$N_r(j)$ and
$N_r(j')$ for some
$r$.  Since a result analogous to Proposition 1.10 was proven in \cite{GGSWW} with $M(j,c,Q)$ and $\pmapq$ replaced by 
$N_r(j)$ and
$N_r(j')$, we will then be able to conclude that $j\simeq j'$.  

Since $\tau$ is both an isometry and a bundle map, it
induces an isometry
$\bt: Q\to Q$.  Moreover, $\tau$
carries the mean curvature vector field for $\pi_j$ to that for $\pi_{j'}$. Thus $\bt_*(\ba)=\ba$, so $\bt$ preserves the level sets of
$\|\ba\|$ in $Q$.  As seen above, each such level set is a union of spheres
$S(t)$ where $t$ ranges over a level set of $q$.  

If $q$ is non-constant, then we may choose a regular value of $q$ so that the corresponding level set of $\|\ba\|$ consists of only
finitely many such spheres.  Choose $t$ in this level set of $q$.  Then $\bt$ carries $S(t)$ to some $S(t')$ with $q(t')=q(t)$.  Moreover $f(t)=f(t')$
since the isometric spheres
$S(t)$ and $S(t')$  must have the same radius.  The restriction of $\tau$ to $\pi_j^{-1}(S(t))$ carries  $\pi_j^{-1}(S(t))$ isometrically to
$\pi_{j'}^{-1}(S(t'))$; i.e., it defines an isometry between $\pn_{f(t)}(j)$ and $\pn_{f(t)}(j')$ in the notation of 1.3(iv).  

Next if $q$ is constant, then $f'$ must also be constant and thus $f$ cannot vanish at both $t_1$ and $t_2$, say $f(t_1)\neq 0$.  In this case, $M(j,c,Q)$ and
$M(j',c,Q)$ are manifolds with boundary; their boundaries are the inverse images under $\pi_j$ and $\pi'_j$, respectively, of either $S(t_1)$ or of $S(t_1)\cup
S(t_2)$, depending on whether $f(t_2)$ is zero. Thus
$\bt$ carries $S(t_1)$ to one of
$S(t_1)$ or
$S(t_2)$.  In either case, the argument in the previous paragraph yields an isometry between $\pn_r(j)$ and $\pn_r(j')$ where $r=f(t_1)$.  

Thus, whether or not $q$ is constant,
we obtain an isometry between 
$\pn_r(j)$ and
$\pn_r(j')$ for some $r$.  As noted above, it follows that $j\simeq j'$.  

Before moving on to the general case, we emphasize some key points in the argument above:  There exists some $t$ (either a regular point of
$q$ or an endpoint of the interval $[t_1,t_2]$) and some $t'$ for which $\tau$ carries  $\pi_j^{-1}(S(t))$ isometrically to
$\pi_{j'}^{-1}(S(t'))$.  Moreover, given such a $t$, then there are only finitely many possibilities for $t'$: either $t'$ is in the same finite level set as
$t$ or else $t$ and $t'$ are both endpoints.  We will also need below the fact that $\bt_*(\ba)=\ba$.

We have reduced the proposition to the following claim:
\smallskip

\noindent {\bf Claim.} Under the hypotheses of the theorem, if $M(j,c,Q)$ and $M(j',c,Q)$ are isometric, then there exists an isometry between them which is
also a bundle map.  

\smallskip

Before proving the claim, we first show that the torus $T$ (see notation 1.3) is a maximal torus in the full isometry group $\Iso(M(j,c,Q))$.  Let $C(T)$ be
the centralizer of
$T$ in 
$\Iso(M(j,c,Q))$  and let
$\tau\in C(T)$.  Since
$\tau$ is an isometry of
$M(j,c,Q)$ which commutes with the action of
$T$, it must also be a bundle map of $M(j,c,Q)$.  Thus by the key  point stated just before the claim,  there exist $t\in[t_1,t_2]$ and $t'$ in a finite
subset of $[t_1,t_2]$
such that
$f(t)=f(t')$ and such that $\tau$ restricts to an isometry
$\sigma$ from 
$\pi_j^{-1}(S(t))$ to $\pi_{j}^{-1}(S(t'))$, with both submanifolds being isometric to $\pn_{f(t)}(j)$.  Moreover, 
$\tau$ is uniquely determined by its restriction $\sigma$ to $\pi_j^{-1}(S(t))$.  Indeed choose a point $p$ in $\pi_j^{-1}(S(t))$.  Then
the horizontal lift $\tilde{A}$ of $\ba$ spans a complement to $T_p(N_{f(t)}(j))$
in
$T_p(M(j,c,Q))$. Since $\bt_*(\ba)=\ba$, we have $\tau_*\tilde{A}=\tilde{A}$.  Thus both the value and the differential of $\tau$ at $p$ are
completely determined by
$\sigma$.  Recalling that an isometry of a connected  manifold is uniquely determined by its value
and differential at a single point, we see that $\sigma$ determines $\tau$.  

Let $r=f(t)$.  The isometry $\sigma$ commutes with the action of $T$ on $N_r(j)$.  Under the hypothesis on $j$ given in the proposition, it was
proven in
\cite{GGSWW} that
$T$ is a maximal torus in the isometry group of
$\pn_r(j)$ for any
$r$, in particular for $r=f(t)$.  This fact together with the fact that there are only finitely many possibilities for $t'$ shows that $T$ is a maximal
connected subset of
$C(T)$; i.e.,
$T$ is a maximal torus in $\Iso(M(j,c,Q))$.  (Aside:  In order to apply the result of \cite{GGSWW}, we need to use the hypothesis on $j$ stated in the
theorem.  Thus we cannot conclude a priori that $T$ is also a maximal torus in $\Iso(M(j',c,Q))$, although we will see below that this is the case.)

 We now prove the
claim.  Suppose
$\tau:M(j,c,Q)\to M(j',c,Q)$ is an isometry.  Then
$\tau$ induces an isomorphism
$\hat{\tau}:\Iso(M(j,c,Q))\to \Iso(M(j',c,Q))$ given by
$\hat{\tau}(\beta)=\tau\beta\tau^{-1}$.  In particular, the maximal tori in  $\Iso(M(j,c,Q))$ and $\Iso(M(j',c,Q))$ must have the same dimension, 
so $T$, viewed as a subgroup of
$\Iso(M(j',c,Q))$, must be a maximal torus.  Since all maximal tori in $\Iso(M(j',c,Q))$ are conjugate in $\Iso(M(j',c,Q))$, we may assume, 
after composing $\tau$ with an isometry of
$M(j',c,Q)$, that $\hat{\tau}$ carries $T$ to $T$.  Equivalently, $\tau$ is a bundle map.  This proves the claim and the proposition follows.

\enddemo

\proclaim{1.12 Proposition  \cite{GW4}}
Let $\dim(\fz)=2$, and let $m=\dim ({\frak v})$ be any positive integer other
than $1,2,3,4$, or $6$. Let $W$ be the real vector space consisting of
all linear maps from $\fz$ to $\so({\frak v})$. Then there is a Zariski open
subset ${\Cal O}$ of $W$ (i.e., ${\Cal O}$ is the complement of the zero locus 
of some non-zero polynomial
function on $W$) such that each
$j \in {\Cal O}$ belongs to a $d$-parameter family of isospectral,
inequivalent elements of $W$. Here $d\geq m(m-1)/2 - [m/2]([m/2]+2)>1$.
In particular, $d$ is of order at least $O(m^2)$.
\endproclaim

\proclaim{Corollary 1.13}   In the notation of 1.3, fix $\fv$, $\fz$, $\fa$, $\Lc$, $c$ and $Q$.  Assume that $\dim(\fz)=2$, that $m=\dim ({\frak v})$ is not
equal to 
$1,2,3,4$, or $6$
 and that $Q$ is either a bounded domain in $R$ or a submanifold of codimension one.  Then for each $j$ in the set ${\Cal O}$ defined in Proposition 1.12, the
manifold $M(j,c,Q)$ lies in a continuous $d$-parameter family of isospectral, non-isometric manifolds, where $d$ is given as in 1.12.  Moreover, if $Q$ is a 
bounded domain, then the boundaries of these isospectral manifolds are also isospectral but not isometric.  
\endproclaim

\demo{Proof}  One of the defining properties of the Zariski set ${\Cal O}$ in 
\cite{GW4} is that the elements $j$ satisfy the hypothesis of Proposition 1.10.  Thus the corollary follows immediately from Propositions 1.10 and 1.12 and
the comments following the statement of Proposition 1.10.
\enddemo

Although the expression for $d$ gives $0$ when $m=6$, explicit examples of
continuous families of isospectral, inequivalent $j$ maps are given in \cite{GW4}, Example 2.4, when $m=6$. These maps also satisfy the hypothesis of
Proposition 1.10 and thus give rise to continuous families of isospectral, non-isometric manifolds $M(j,c,Q)$.

\medskip
\heading Locally homogeneous examples\endheading

Throughout this section, we use the notation of 1.3, with $Q$ chosen to be a bounded domain in $R$, for example a geodesic ball.  Recall that $M(j,c,Q)$ is
locally homogeneous and its local geometry is that of $G(j,c)$.

The following proposition is a special case of the more general classification of homogeneous manifolds of non-positive curvature, carried out in \cite{Hz}
for strictly negative curvature and in \cite{AW1,2} for nonpositive curvature.

\proclaim{2.1 Proposition}(\cite{Hz} or \cite{AW2}, Proposition 8.5)  In the notation of 1.3, for each choice of $j$, there exists a constant
$\lambda(j)$ such that $G(j,c)$ has strictly negative curvature when $c>\lambda(j)$, non-positive curvature when $c=\lambda(j)$ and mixed curvature when
$c<\lambda(j)$.  Moreover, as
$c$ approaches $\infty$, the maximum of the sectional curvature of $G(j,c)$ approaches $-\infty$.
\endproclaim

\proclaim{2.2 Theorem}  In every dimension $n\geq 8$, there exist continuous isospectral deformations of locally homogeneous negatively curved manifolds 
with boundary which are not locally isometric. 
Moreover, given any constant
$\kappa >0$, we can choose the isospectral metrics so that their curvature is bounded above by $-\kappa$.
\endproclaim

\demo{Proof}  Let $k=2$, let $m\geq 5$ and let $\{j_t\}$ be a family of isospectral, inequivalent maps $j_t:\fz\to\so(\fv)$ as in Proposition 1.12 or as in
the comments following Corollary 1.13. 
 Then for any choice of $c>0$, the
manifolds
$(M(j_t,c,Q))$ are isospectral but not locally isometric. Proposition 2.1 allows us to adjust the constant $c$ to achieve any desired curvature bound.
\enddemo
 
\remark{Remark 2.3}  If the constant $\lambda(j_t)$ depends non-trivially on $t$, then we can choose $c$ so that $c>\lambda(j_t)$ for some choices of $t$ and
$c<\lambda(t)$ for some other choices of $t$.  In this case, we obtain isospectral deformations in which some of the metrics have negative curvature and
others mixed.  Unfortunately, due to the complexity of the curvature expressions, we have not been able to compute the constants $\lambda(j_t)$ for any of
the deformations.  It is very likely, however, that for generic deformations, $\lambda(j_t)$ does vary with $t$.  The only situation in which we have
computed the constants thus far is for the pair of isospectral $j$ maps in Example 2.4 below.  In this case, the constants are indeed different.  This pair of
maps does not belong to a continuous family of isospectral,  inequivalent maps, however.
\endremark

\definition{2.4 Example}  A two-step nilmanifold $H(j)$ defined as in 1.3(ii) is said to be of {\it Heisenberg type}, as defined by A. Kaplan \cite{K}, if
$j^2=-Id$. For simplicity, we will say that
$j$ is of Heisenberg type in this case.  Damek and
Ricci \cite{DR} proved that if $j$ is of Heisenberg type, then $G(j,1)$ is a harmonic manifold.  Included among such spaces are all the rank one symmetric
spaces of non-compact type as well as the first known examples of non-symmetric harmonic manifolds.  

Observe that if
$j,j':\fz\to\so(\fv)$ are both of Heisenberg type, then they are necessarily isospectral.  We now consider a specific example.  Let
$\frak z$ be the purely imaginary quaternions with the
standard inner product, and let $\frak v$
be the orthogonal direct sum of $l$ copies of the quaternions, viewed as a
$4l$-dimensional real vector space with the standard inner product.  Choose
non-negative integers $a$ and $b$
with $l=a+b$.  Define the map
$j_{a,b}:{\frak z}
\to\text{so}({\frak v})$ by    
$$j_{a,b}(p)(q_1,\dots, q_a, q'_1,\dots,q'_b)=(pq_1,\dots,pq_a,q'_1 p,\dots,q'_b
p)$$ where $pq_i$ and
$q'_j p$ denote quaternionic multiplication.  The maps $j_{(a,b)}$ and $j_{(a',b')}$ are both of Heisenberg type, and they are equivalent
if and only if $(a',b')$ is a permutation of $(a,b)$. Thus when $(a',b')$ is not a permutation of $(a,b)$ but $a+b=a'+b'$, the manifolds $M(j_{(a,b)},c,Q)$ and
$M(j_{(a',b')},c,Q)$ are isospectral but not isometric for any choice of $c$ and bounded domain $Q$.  In particular, choosing $c=1$ and choosing any pair
$(a,b)$ of positive integers, we obtain isospectral manifolds
$M(j_{(a+b,0)},1,Q)$ and
$M(j_{(a,b)},1,Q)$, the first of which is a domain with boundary in the symmetric space $G(j_{(a+b,0)},1)$ (quaternionic hyperbolic space) and thus has
parallel curvature tensor while the second has non-parallel curvature.

We now show that $\lambda(j_{(2,0)})\neq \lambda(j_{(1,1)})$.  The quaternionic hyperbolic space
$G(j_{(2,0)},1)$ has curvature bounded above by -1/4.  Thus
$\lambda(j_{(2,0)})>1$. Damek \cite{D} showed that if $j$ is of Heisenberg type, then $G(j,1)$ always has non-positive curvature, and it has some two planes
of zero curvature if and only if there exist vectors $X,Y\in\fv$ such that $[X,Y]=0$ and $j(\fz)(X)\cap j(\fz)(Y)$ is non-empty.  Consider $j=j_{(1,1)}$.  In
the Lie algebra $\fh(j)$, the subspace $W$ of $\fv$ given by $\{(q,q):q\in Q\}$ is an abelian subalgebra of $\fh(j)$.  Letting $q$ and $q'$ be distinct
non-zero purely imaginary quaternions, then the elements $X=(q,q)$ and $Y=(q',q')$ of
$\fv$ satisfy Damek's condition for the existence of zero curvature.   Thus $G(j_{(1,1)},1)$ has some zero curvature, and so $\lambda(j_{(1,1)})=1$. 
Choosing
$c$ so that
$\lambda(j_{(2,0)})>c>1$, we find that $M(j_{(2,0)},c,Q)$ has strictly negative curvature while $M(j_{(1,1)},c,Q)$ has mixed curvature.
\enddefinition

We next consider conditions for constant Ricci curvature and for constant scalar curvature.

\proclaim{2.5 Proposition}\cite{EH} We use the notation of 1.3.  The solvable Lie group $G(j,c)$ with its left-invariant metric is
Einstein if and only  if $c=1$ and both the following conditions are satisfied:

(i) The map $j:\fz\to j(\fz)\subset\so(\fv)$ is a linear isometry relative to the Riemannian inner
product $\langle\,,\,\rangle$ on $\fz$ and the inner product $(\,,\,)$ on $\so(\fv)$ given by
$(\alpha,\beta)=-\frac{1}{m}\tr(\al\beta)$, where $m=\dim(\fv)$.

(ii)  Letting $\{Z_1,\dots ,Z_k\}$ be an orthonormal basis of $\fz$, then $\sum_{i=1}^k
j(Z_i)^2$ is a scalar operator. (This condition is independent of the choice of orthonormal basis.)
\endproclaim

We will now construct a pair  of isospectral maps $j,j'$ such that $G(j,1)$ is Einstein but $G(j',1)$ is not.  For any choice of bounded domain $Q$, this will
then give us a pair of isospectral manifolds $M(j,1,Q)$ and
$M(j',1,Q)$ with boundary such that $M(j,1,Q)$ has constant Ricci curvature, but $M(j',1,Q)$ does not.

\definition{2.6 Example}  We take $\fz=\R^3$ with the inner product 
$\langle Z,W\rangle=\frac{2}{3}(Z,W)$, where
$(\,,\,)$ is the standard inner product.  We take $\fv=\R^6$ with its standard inner product.  To define $j$, view $\fv$ as $\R^3\times\R^3$, and for
$Z\in\fz=\R^3$, and
$(U,V)\in\fv$, set 
$$J(Z)(U,V)=(Z\times U,Z\times V)$$
where $\times$ denotes the cross product in $\R^3$.  Then the eigenvalues of $j(Z)$ are $\pm\sqrt{\frac{3}{2}}\|Z\|\sqrt{-1}$ and $0$, each occurring with multiplicity 2.  It is
straightforward to verify that $j$ satisfies the conditions of Proposition 2.5, so that $G(j,1)$ is Einstein and $M(j,1,Q)$ has constant Ricci curvature.

To define $j'$, view $\fv=\R^6$  as $Q\times \R^2$ where $Q$ denotes the quaternions, and view $\fz=\R^3$ as the space of purely imaginary quaternions.  For $Z\in\fz$ and $(U,V)\in\fv$
with $U\in Q$ and $V\in\R^2$, set 
$$j'(Z)(U,V)=(ZU,0)$$
where $ZU$ denotes quaternionic multiplication.  Then the eigenvalues of $j'(Z)$ are also $\pm\sqrt{\frac{3}{2}}\|Z\|\sqrt{-1}$ and $0$, each occurring with multiplicity 2.  However,
now all the $j'(Z)$, as $Z$ varies over $\fz$, have the same zero eigenspace, so the second condition in Proposition 2.5 is not satisfied.  Thus $G(j',1)$ is
not Einstein, and
$M(j',1)$ does not have constant Ricci curvature.  
\enddefinition

\proclaim{2.7 Theorem}  In the notation of 1.3, the closed manifold $ N(j)$ has constant scalar curvature if and only if $j$ satisfies condition
(ii) of Proposition 2.5.
\endproclaim

\demo{Proof}  The Lie group exponential map $\exp:\fh(j)\to H(j)$ is a diffeomorphism.  Define global coordinates $H(j)\to\fv +\fz$ on $H(j)$ by $\exp(x+z)\to x+z$
for $x\in\fv$ and
$z\in\fz$.  We then obtain a diffeomorphism of $\bhj=\Lc\bs H(j)$ with $\fv\times T$ by recalling that $T\simeq \Lc\bs\fz$.  We use this diffeomorphism to parametrize $\bhj$,
denoting points in $\bhj$ as $(x,\bar{z})$ with $x\in\fv$ and $\bar{z}\in T$.    

Since $\bhj$ is homogeneous, it has constant scalar curvature $\tau$.  In
\cite{GGSWW}, the scalar curvature of the submanifold $N(j)$ is computed:
$$scal(x,\bar{z})=\tau(x,\bar{z})+(m-1)(m-2)-\sum_{i=1}^k j(Z_i)^2(x,x)$$
where $\{Z_1,\dots ,Z_k\}$ is an orthonormal basis of
$\fz$.  
The theorem now follows.
\enddemo

\definition{2.8 Example}  In Example 2.6, we constructed a pair of isospectral maps $j,j':\R^3\to\so(6,\R)$ such that $j$ satisfies condition (ii) of Proposition
2.5 but
$j'$ does not.  These maps give rise to 8-dimensional closed manifolds $ N(j)$ and $ N(j')$, such that the first has constant scalar curvature and
the second does not.
\enddefinition

\medskip
\heading Locally inhomogeneous examples\endheading

In this section we study the geometry of the manifolds of the form $\pmaq$ in the notation 1.3 for which $Q$ is a submanifold of $R$ of codimension 
one. For simplicity we will focus on the special class of such submanifolds specified in 3.2 below.
Before considering the manifolds $\pmaq$,   we first describe the connection and curvature on the ambient space
$\bgj$.  Recall that
$\bgj$ is locally homogeneous and locally isometric to $G(j,c)$, so it suffices for this purpose to consider $G(j,c)$.

\proclaim{Proposition 3.1}  We use the notation of 1.3.  For $X$, $X^*\in\fv$ and $Z,Z^*\in\fz$, we have:
\item{(a)} $\nabla_XX^*=\frac{1}{2}[X,X^*]+\frac{c^2}{2}\langle X,X^*\rangle A.$
\item{(b)} $\nabla_ZZ^*=c^2\langle Z,Z^*\rangle A.$
\item{(c)} $\nabla_XZ=\nabla_ZX=-\frac{1}{2}J_Z(X).$
\item{(d)} $\nabla_XA=-\frac{1}{2}X.$
\item{(e)} $\nabla_ZA=-Z.$
\item{(f)} $\nabla_A =0.$
\endproclaim

The proposition is a straightforward computation using the fact that for left-invariant vector fields, we have
$$2\langle \nabla_UV,W\rangle =\langle [U,V],W\rangle +\langle [W,U],V\rangle +\langle [W,V],U\rangle $$
and recalling that $\|A\|=\frac{1}{c}$.

\definition{3.2.  A special class of submanifolds}  Let $Q$ be a submanifold of $R$ of codimension one defined as in 1.3(iv).  As in 
1.11, $Q$ may be expressed  as a level set of a function
$\phi(t,X)=\|X\|-f(t)$.  For simplicity, we now restrict attention to the case that the function $f$ is constant, say $f\equiv r$ with $r>0$.  (In particular, $Q$ has non-trivial
boundary.)  Since
$Q$ is uniquely defined by the constant $r$, we will write $\pmar$ for $\pmaq$ in this case.
\enddefinition

\definition{3.3.  The Weingarten map}  Let $Q$ be given as in 3.2.  A unit normal vector field
$\n_Q$ along
$Q$ is given by 
$\n_Q(t,X)=\frac{X}{r}$ in the notation of 1.11(a).  Since $\pmar=\pi_j^{-1}(Q)$, a unit normal vector field $\n$ along $\pmar$ in $\bgj$ is horizontal with respect to the submersion
$\pi_j$ and is given as follows:  For $a\in\pmar$ with $\pi_j(a)=(t,X)$, we have
$\n_a=\frac{X_a}{r}$ where $X\in\fv$ is viewed as an invariant vector field on $\bgj$.

For $a\in\pmar$ as above, the tangent space to $\pmar$ at $a$ is spanned by $A_a$ together with $\{Z_a:Z\in\fz\}$ and $\{Y_a:Y\in\fv, Y\perp X\}$, where $U_a$
denotes the value at $a$ of an invariant vector field $U$ on $\bgj$.  For simplicity, we will drop the subscripts $a$. 
The Weingarten map $B$ of $\pmar$, defined by $B(U)=\nabla_U(\n)$ for $U$ a tangent vector to $\pmar$, is given at $a$ by:
$$B(Y)=\frac{t^\ot}{r}Y+ \frac{1}{2r}[Y,X],$$
$$B(Z)=-\frac{1}{2r}J_Z(X),$$
$$B(A)=0$$
for $Y\in\fv$ with $Y\perp X$ and for $Z\in\fz$.  
\enddefinition

The key point here is that the Weingarten map does not depend on the parameter $c$.  This fact can be seen even without the explicit formula in 3.3 
by the following observations: (i) when $X\in\fv$
and $U\in\fg(j)$ is perpendicular to $X$, then the invariant vector field $\nabla_UX$ is independent of $c$ (as can be seen from Proposition 3.1) and (ii) the normal vector field $\n$
to $\pmar$ takes all its values in $\fv$.

\proclaim{Theorem 3.4}  We use the notation of 3.2 and 1.3.  Given $j$ as in 1.3(i) and r as in 3.2, there exists a constant $\lambda(j,r)>0$ such that $\pmar$ has strictly negative
curvature when $c>\lambda(j,r)$.  Moreover as $c$ approaches $\infty$, the maximum of the sectional curvature of $\pmar$ approaches $-\infty$.
\endproclaim

\demo{Proof}  Let $R$, respectively $\tilde{R}$, denote the curvature tensor of $\bgj$, respectively $\pmar$.  Then for $X,Y$ tangent vectors to $\pmar$ at a point $a\in\pmar$, we have
$$\langle\tilde{R}(X,Y)Y,X\rangle=\langle R(X,Y)Y,X\rangle +\langle B(X),X\rangle\langle B(Y),Y\rangle -\langle B(X),Y\rangle^2.$$
Since the Weingarten map $B$ is independent of $c$, the theorem follows immediately from Proposition 2.1.
\enddemo

\proclaim{Corollary 3.5}  Suppose $Q$ is a submanifold of $R$ defined as a level set of a function
$\phi(t,X)=\|X\|-f(t)$.  If $f'$ and $f''$ are bounded sufficiently close to zero, then for each $j$, there exists a constant $\lambda(j)>0$ such that $\pmaq$
has strictly negative curvature when $c>\lambda(j)$.  Moreover as $c$ approaches $\infty$, the maximum of the sectional curvature of $\pmaq$ approaches
$-\infty$.
\endproclaim

The corollary follows from the theorem by a continuity argument.  

\remark{3.6 Remark}  One would not expect to obtain a similar result when $f'$ is allowed to vary too greatly.  Indeed, in the extreme case when $f$ is zero on
the endpoints of its domain
$[t_1,t_2]$ and has infinite slope at these points, then $\pmaq$ is a closed manifold which admits a non-trivial isometric action by a torus.  Since Bochner's
Theorem states that a closed manifold of negative Ricci curvature cannot admit a nontrivial Killing field, $\pmaq$ cannot have negative Ricci curvature, let
alone negative sectional curvature, in this case.
\endremark

\proclaim{Proposition 3.7}  In the notation of 3.2, the Riemannian manifold $\pmar$ has non-constant scalar curvature and thus is locally inhomogeneous.
\endproclaim

\demo{Proof}  Letting $\rho$ and $\Ric$ denote the scalar curvature and Ricci tensor on $\bgj$ and $\tilde{\rho}$ the scalar curvature on $\pmar$, then
$$\tilde{\rho}=\rho-2\Ric(\n,\n)+(\tr(B))^2-\tr(B^2).$$  
Since $\bgj$ is locally homogeneous, $\rho$ is constant.  Using coordinates $(t,X)$ on $R$ as in 1.11(a), we see from 3.3 that as $a$ varies over $\pmar$, the normal vector field $\n$
depends only on the $X$ coordinate of $\pi(a)$ and thus $\Ric(\n,\n)$ is independent of the coordinate $t$.  On the other hand, 3.3 also shows that
$(\tr(B))^2-\tr(B^2)$ depends non-trivially on $t$.  Thus $\tilde{\rho}(a)$ depends non-trivially on the $t$ coordinate of $\pi(a)$.
\enddemo

\enddocument